\documentstyle[12pt]{article}
\raggedright
\begin{document}
\setlength{\baselineskip}{24pt}
\setlength{\parindent}{.5in}
\setlength{\footskip}{.3in}
\begin{center}
{\Large Topological Partition Relations of the Form $\omega^{*}\rightarrow(Y)^{1}_{2}$}\\
\end{center}
\medskip
\begin{center}
{\large W. W. Comfort}~$^{a,b,c}$, Wesleyan University (USA)\\
{\large Akio Kato}, National Defense Academy (Japan)\\
{\large Saharon Shelah}~$^{d}$, Hebrew University (Israel)\\
\end{center}
\begin{center}
ABSTRACT\\
\end{center}

Theorem. The topological partition relation $\omega ^{*}\rightarrow(Y)^{1}_{2}$

(a) fails for every space $Y$ with $|Y|\geq 2^{\rm \bf c}$;

(b) holds for $Y$ discrete if and only if $|Y|\leq$ {\bf c};

(c) holds for certain non-discrete $P$-spaces $Y$;

(d) fails  for $Y=\omega\cup\{p\}$ with $p\in\omega^{*}$;

(e) fails for $Y$ infinite and countably compact.\\
\vspace{12pt}

AMS Classification Numbers:   Primary 54A25, 54D40, 05A17; Secondary 54A10,
54D30,
04A20

\vspace{.6in}
\begin{center}
[Footnotes]
\end{center}

$^{a}$ Presented at the Madison Conference by this co-author.

$^{b}$ Member, New York Academy of Sciences.

$^{c}$ This author gratefully acknowledges support received from the Technische
Hochschule Darmstadt and from the Deutscher Akademischer Austauschdienst
(=DAAD) of the Federal Republic of Germany.

$^{d}$ This author acknowledges partial support received from the Basic
Research Fund of the Israeli Academy of Sciences, Publ. 440.
\newpage

\S1. Introduction.

For topological space $X$ and $Y$ we write $X\approx Y$ if $X$
and $Y$ are homeomorphic, and we write $f:X\approx Y$ if $f$ is a
homeomorphism of $X$ onto $Y$. The ``topological inclusion relation"
is denoted by $\subseteq_{h}$; that is, we write $Y\subseteq_{h} X$
if there is $Y'\subseteq Y$ such that $Y\approx Y'$.

The symbol $\omega$ denotes both the least infinite cardinal and the
countably infinite discrete space; the Stone-\v{C}ech remainder
$\beta(\omega)\backslash\omega$ is denoted $\omega^{*}$.

For a space $X$ we denote by $wX$ and $dX$ the weight and density
character of $X$, respectively. Following [7], for
$A\subseteq\omega$ we write $A^{*}=({\rm
cl}_{\beta(\omega)}A)\backslash\omega$.

For proofs of the following statements, and for other basic information
on topological and combinatorial properties of the space $\omega^{*}$,
see [7], [3], [12].

1.1. Theorem. (a) \{cl$_{\beta(\omega)}A:A\subseteq\omega$\} is a
basis for the open sets of $\beta(\omega)$; thus
$w(\beta(\omega))=\bf c$.

(b) There is an (almost disjoint) family $\cal A$ of subsets of $\omega$
such that $|\cal A|=\bf c$ and $\{A^{*}:A\in\cal A\}$ is pairwise
disjoint.

(c) $\omega^{*}$ contains a family of $2^{\bf c}$-many pairwise disjoint
copies of $\beta(\omega)$.

(d) Every infinite, closed subspace $Y$ of $\omega^{*}$ contains a copy
of $\beta(\omega)$, so $|Y|=|\beta(\omega)|=2^{\bf c}$.~~$\Box$

For cardinals $\kappa$ and $\lambda$ and topological spaces $X$
and $Y$, the symbol $X\rightarrow(Y)^{\kappa}_{\lambda}$ means that
if the set $[X]^{\kappa}$ of all $\kappa$-membered subsets of $X$
is written in the form $[X]^{\kappa}=\cup_{i<\lambda}P_{i}$, then
there are $i<\lambda$ and $Y'\subseteq X$ such that $Y\approx
Y'$ and $[Y']^{\kappa}\subseteq P_{i}$. Our present primary interest is in
topological arrow relations of the form $X\rightarrow(Y)^{1}_{2}$ (with $X=\omega^{*}$). For
spaces $X$ and $Y$, the relation $X\rightarrow(Y)^{1}_{2}$ reduces to this: if
$X=P_{0}\cup P_{1}$, then either $Y\subseteq_{h}P_{0}$ or $Y\subseteq_{h}P_{1}$.

The relation $X\rightarrow(Y_{0},Y_{1})^{1}$ indicates that if
$X=P_{0}\cup P_{1}$, then either $Y_{0}\subseteq_{h}P_{0}$ or
$Y_{1}\subseteq_{h}P_{1}$.

It is obvious that if $X$ and $Y$ are spaces such that $Y\subseteq_{h}X$
fails,
then $X\rightarrow(Y)^{1}_{2}$ fails.

By way of introduction it is enough here to observe that the classical
theorem of F.~Bernstein, according to which there is a subset $S$ of the real
line {\bf R} such that neither $S$ nor its complement
${\bf R}\backslash S$ contains an uncountable closed set,
is captured by the assertion that the relation
{\bf R}~$\rightarrow(\{0,1\}^{\omega})^{1}_{2}$ fails; in the positive
direction, it is easy to see that the relation
{\bf Q}~$\rightarrow({\bf Q})^{1}_{2}$ holds for
{\bf Q} the space of rationals.

For a report on the present-day ``state of the art" concerning
topological partition relations, and for references to the literature and
open questions, the reader may consult [14], [15], [16].

This paper is organized as follows. \S2 shows that $\omega^{*}\rightarrow(Y)^{1}_{2}$
fails for every infinite compact space $Y$. \S3 characterizes those
discrete spaces $Y$ for which $\omega^{*}\rightarrow(Y)^{1}_{2}$,
and \S4 shows that $\omega^{*}\rightarrow(Y)^{1}_{2}$ holds for
certain non-discrete spaces $Y$. \S5 shows that
$\omega^{*}\rightarrow(Y)^{1}_{2}$ fails for spaces of the
form $Y=\omega\cup\{p\}$ with $p\in\omega^{*}$, hence fails for every
infinite countably compact space $Y$. The results of \S\S2--5 prompt
several questions, and these are given in \S6.

We are grateful to Jan van Mill, K. P. S. Bhaskara Rao, and W. A. R.
Weiss for helpful conversations.

We announced some of our results in the abstract [2]. See also [1] for
related results.

\S2. $\omega^{*}\not\rightarrow(Y)^{1}_{2}$ for $|Y|\geq2^{\bf c}$.

2.1. Lemma. If $Y\subseteq_{h}\omega^{*}$ then $|\{A\subseteq\omega^{*}:A\approx
Y\}|=2^{\bf c}$.

Proof. The inequality $\geq $ is immediate from Theorem 1.1(c). For
$\leq $, it is enough to fix (a copy of) $Y\subset\omega^{*}$ and to
notice that since $dY\leq wY\leq w(\omega^{*})=\bf c$ (by Theorem
1.1(a)), the number of continuous functions from $Y$ into $\omega^{*}$
does not exceed $|(\omega^{*})^{dY}|\leq(2^{\bf c})^{\bf c}=2^{\bf
c}$.~~$\Box$

2.2. Theorem. If $Y$ is a space such that $|Y|\geq2^{\bf c}$, then
$\omega^{*}\not\rightarrow(Y)^{1}_{2}$.

Proof. We assume $Y\subseteq_{h}\omega^{*}$ (in particular we assume
$|Y|=|\omega^{*}|=2^{\bf c}$) since otherwise $\omega^{*}\not\rightarrow(Y)^{1}_{2}$ is
obvious. Following
Lemma 2.1 let $\{A_{\xi}:\xi<2^{\bf c}\}$ enumerate
$\{A\subseteq\omega^{*}:A\approx Y\}$, choose distinct $p_{0}, q_{0}\in
A_{0}$ and recursively, if $\xi<2^{\bf c}$ and $p_{\eta}, q_{\eta}$
have been chosen for all $\eta<\xi$ choose distinct

$p_{\xi},
q_{\xi}\in A_{\xi}\backslash(\{p_{\eta}:\eta<\xi\}\cup\{q_{\eta}:\eta<\xi\})$.\\

\noindent It is then clear, writing

$P_{0}=\{p_{\xi}:\xi<2^{\bf c}\}$ and $P_{1}=\omega^{*}\backslash
P_{0}$,\\
\noindent that the relations $Y\subseteq_{h}P_{0}$ and $Y\subseteq_{h}P_{1}$ both
fail.~~$\Box$

The following statement is an immediate consequence of Theorems 2.2 and
1.1(d).

2.3. Corollary. The relation $\omega^{*}\rightarrow(Y)^{1}_{2}$ fails for every infinite compact space
$Y$.~~$\Box$

By less elementary methods we strengthen Corollary 2.3 in Theorem 5.14 below. 

\S3. Concerning the Relation $\omega^{*}\rightarrow(Y)^{1}_{2}$ for $Y$ Discrete.

The very simple result of this section, included in the interest of
completeness, shows for discrete spaces $Y$ that $\omega^{*}\rightarrow(Y)^{1}_{2}$
if and only if $Y\subseteq_{h}\omega^{*}$.

3.1. Theorem. For a discrete space $Y$, the following conditions are
equivalent.

(a) $|Y|\leq{\bf c}$;

(b) $\omega^{*}\rightarrow(Y)^{1}_{\bf c}$;

(c) $\omega^{*}\rightarrow(Y)^{1}_{2}$;

(d) $Y\subseteq_{h}\omega^{*}$.

Proof. (a) $\Rightarrow$ (b).
[Here we profit from a suggestion offered by the referee.] Given
$\omega^{*}=\cup_{i<{\bf c}}~P_{i}$, recall from [10](2.2) or [12](3.3.2)
this theorem of Kunen: there is a matrix
$\{A^{\xi}_{i}:\xi<{\bf c},~i<{\bf c}\}$
of clopen subsets of $\omega^{*}$ such that

(i) for each $i<{\bf c}$ the family $\{A^{\xi}_{i}:\xi<{\bf c}\}$
is pairwise disjoint, and

(ii) each $f\in{\bf c}^{\bf c}$ satisfies
$\cap_{i<{\bf c}}~A^{f(i)}_{i}\neq\emptyset$.\\
\noindent Now if one of the sets $P_{i}$ meets $A^{\xi}_{i}$
for each $\xi<{\bf c}$
(say $p_{\xi}\in A^{\xi}_{i}$) then the discrete set
$D=\{p_{\xi}:\xi<{\bf c}\}$ satisfies
$Y\subseteq_{h} D\subseteq P_{i}$; otherwise for each $i<{\bf c}$
there is $f(i)$ such that
$P_{i}\cap A^{f(i)}_{i}=\emptyset$, so
$\emptyset\neq\cap_{i<{\bf c}}~A^{f(i)}_{i}\subseteq\omega^{*}\backslash\cup_{i<{\bf c}}~P_{i}$.

That (b) $\Rightarrow$ (c) and (c) $\Rightarrow$ (d)
and clear.

(d) $\Rightarrow$ (a). Theorem 1.1(a) gives $|Y|=wY\leq
w(\beta(\omega))={\bf c}$.

\S4. $\omega^{*}\rightarrow(Y)^{1}_{2}$ for Certain Non-Discrete $Y$.

For an infinite cardinal $\kappa$ we denote by $P_{\kappa}$ the
ordinal space $\kappa+1=\kappa\cup\{\kappa\}$ topologized to be
``discrete below $\kappa$"
and with a neighborhood base at $\kappa$ the same as in the usual
interval topology. That is, a subset $U$ of $\kappa+1$ is open in
$P_{\kappa}$ if and only if either $U\subseteq\kappa$ or some
$\xi<\kappa$ satisfies $(\xi,\kappa]\subseteq U$.

4.1. Theorem. For cardinals $\kappa\geq\omega$
and $m_{0},m_{1}<\omega$, the space $P_{\kappa}$ satisfies
$P_{\kappa}^{m_{0}+m_{1}}\rightarrow(P_{\kappa}^{m_{0}},P_{\kappa}^{m_{1}})^{1}$.

Proof. Let $P^{I}=X_{0}\cup X_{1}$ and $|I|=m_{0}+m_{1}$ and suppose
without loss of generality that the point $c=\langle c_{i}\rangle_{i\in I}$
with $c_{i}=\kappa$ (all $i\in I)$ satisfies $c\in X_{0}$. Let $I=I_{0}\cup I_{1}$
with $|I_{0}|=m_{0}$, $|I_{1}|=m_{1}$,
and set $D=P_{\kappa}\backslash\{\kappa\}$, and for $x\in D^{I_{0}}$ define

$S(x)=\{x\}\times\{y\in P_{\kappa}^{I_{1}}:{\rm max}\{x_{i}:i\in I_{0}\}
< {\rm min}\{y_{i}:i\in I_{1}\}\}$.\\
\noindent If some $x\in D^{I_{0}}$ satisfies $S(x)\subseteq X_{1}$ we have
$P_{\kappa}^{m_{1}}\approx S(x)\subseteq X_{1}$
and the proof is complete. Otherwise for each
$x\in D^{I_{0}}$ there is $p(x)\in S(x)\cap X_{0}$ and then

$P_{\kappa}^{m_{0}}\approx\{p(x):x\in D^{I_{0}}\}\cup\{c\}\subseteq X_{0}$,\\
\noindent as required.~~$\Box$

4.2. Corollary. Every infinite cardinal $\kappa$ satisfies
$P_{\kappa}\times P_{\kappa}\rightarrow (P_{\kappa})^{1}_{2}.$~~$\Box$

We say as usual that a topological space
$X=\langle X,{\cal T}\rangle$ is a $P$-space if each
${\cal U}\subseteq{\cal T}$
with $|{\cal U}|\leq\omega$ satisfies $\cap{\cal U}\in{\cal T}$, Since (clearly)
$P_{\kappa}$ is a non-discrete $P$-space if and only if cf$(\kappa)>\omega$,
the following theorem shows the existence of a nondiscrete $Y$ such that
$X\rightarrow(Y)^{1}_{2}$.

4.3. Theorem. Let $\omega_{1}\leq\kappa\leq{\bf c}$ satisfy
cf$(\kappa)>\omega$. Then $\omega^{*}\rightarrow(P_{\kappa})^{1}_{2}$.

Proof. It is a theorem of E. K. van Douwen that every $P$-space $X$
such that $wX\leq{\bf c}$ satisfies $X\subseteq_{h}\omega^{*}$. (For a proof of
this result see [4] or [12].) Thus for $\kappa$ as hypothesized we
have $P_{\kappa}\times P_{\kappa}\subseteq_{h}\omega^{*}$, so the relation
$\omega^{*}\rightarrow(P_{\kappa})^{1}_{2}$ is immediate from Corollary
4.2.~~$\Box$

4.4. Remarks.
(a) The following simple result, suggested by the proof of Theorem 4.2,
is peripheral to the principal thrust of our paper. Here as usual for a space
$X=\langle X,\cal T\rangle$ we denote by $PX=\langle PX,P\cal T\rangle$
the set $X$ with the smallest topology $P\cal T$ such that
$P\cal T\supseteq\cal T$ and $PX$ is a $P$-space; thus
$\{\cap\cal U:\cal U\subseteq\cal T, |\cal U|\leq\omega\}$
is a base for $P\cal T$.

Theorem. For a $P$-space $Y$, the following conditions are equivalent.

(i) $\omega^{*}\rightarrow(Y)^{1}_{2}$;

(ii) $\{0,1\}^{\bf c}\rightarrow(Y)^{1}_{2}$;

(iii) $P(\omega^{*})\rightarrow(Y)^{1}_{2}$;

(iv) $P(\{0,1\}^{\bf c})\rightarrow(Y)^{1}_{2}$.

Proof. The implications (iii) $\Rightarrow$ (iv) $\Rightarrow$ (i)
$\Rightarrow$(ii) follow respectively from the inclusions
$P(\omega^{*})\subseteq_{h}P(\{0,1\}^{\bf c})\subseteq_{h}\omega^{*}\subseteq_{h}
\{0,1\}^{\bf c}$. (Of these three inclusions
the third follows from Theorem 1.1, the first from the third, and the second
from van Douwen's theorem cited above.) That (ii) $\Rightarrow$ (iii)
follows from $P(\{0,1\}^{\bf c})\subseteq_{h}\omega^{*}$
(whence $P(\{0,1\}^{\bf c})\subseteq_{h} P(\omega^{*}))$ and the case
$A=\{0,1\}^{\bf c}$, $B=Y=PY$ of this general observation:
if $A\rightarrow(B)^{1}_{2}$ then $PA\rightarrow(PB)^{1}_{2}$.~~$\Box$

(b) We note in passing the following result, from which (with 4.1) it
follows that for $\kappa\geq\omega$ the space $P_{\kappa}$ satisfies
$P_{\kappa}^{2^{n}}\rightarrow(P_{\kappa})^{1}_{n+1}$.

Theorem. Let $S$ be a space such that
$S^{m_{0}+m_{1}}\rightarrow(S^{m_{0}},S^{m_{1}})^{1}$ for $m_{0},m_{1}<\omega$.
Then $S^{2^{n}}\rightarrow(S)^{1}_{n+1}$ for $n<\omega$.\hfill{(*)}\\

Proof. Statement 
(*) is trivial when $n=0$, and is given by the case $m_{0}=m_{1}=1$ of the hypothesis when
$n=1$.

Now suppose (*) holds for $n=k$, and let
$S^{2^{k+1}}=\cap_{i=0}^{k+1}~X_{i}$. With
$Y_{0}=X_{0}$ and $Y_{1}=\cup_{i=1}^{k+1}~X_{i}$, it follows from
$S^{2^{k}+2^{k}}\rightarrow(S^{2^{k}},S^{2^{k}})$
that there is $T\subseteq S^{2^{k+1}}$ such that
$T\approx S^{2^{k}}$ and either
$T\subseteq Y_{0}$ or $T\subseteq Y_{1}$. In the first case we have
$S\subseteq_{h} T\subseteq X_{0}$, and in the second case from
$T\subseteq\cup_{i=1}^{k+1}~X_{i}$ and (*) at $k$ there exists $i$
such that $1\leq i\leq k+1$ and $S\subseteq_{h} X_{i}$, as required.~~$\Box$

(c) The method of proof of 4.1
and 4.2 applies to many spaces other
than those of the form $P_{\kappa}$. The reader
may easily verify, for example,
denoting by $C_{\kappa}$ the one-point compactification of the
discrete space $\kappa$, that
$C_{\kappa}\times C_{\kappa}\rightarrow(C_{\kappa})^{1}_{2}$,
and hence $\{0,1\}^{\kappa}\rightarrow(C_{\kappa})^{1}_{2}$, for all
$\kappa\geq\omega$. For a proof due to S.~Todor\v{c}evi\'{c}
of a much stronger topological partition relation, namely
$\{0,1\}^{\kappa}\rightarrow(C_{\kappa})^{1}_{{\rm cf}(\kappa)}$,
see Weiss~[15].

\S5. $\omega^{*}\not\rightarrow(Y)^{1}_{2}$ for $Y$ Infinite and Countably Compact.

To prove this result, we show first that the relation
$\omega^{*}\rightarrow(\omega\cup\{p\})^{1}_{2}$ fails for every $p\in\omega^{*}$. While this can be proved directly by
combinatorial arguments, we find it convenient (given $p\in\omega^{*}$)
to introduce and use as a tool a new topology ${\cal T}(p)$ on
$\omega^{*}$.

Given $f:\omega\rightarrow\omega^{*}$, we denote by
$\overline{f}:\beta(\omega)\rightarrow\omega^{*}$ the Stone extension of
$f$. For $X\subseteq\omega^{*}$ we set

$X^{p}=X\cup\{\overline{f}(p):f:\omega\approx f[\omega]\subseteq
X\}$;\\
\noindent that is, $X^{p}$ is $X$ together with its ``$p$-limits
through discrete countable sets."

5.1. Lemma. There is a topology ${\cal T}(p)$ for $\omega^{*}$ such
that each $X\subseteq\omega^{*}$ satisfies: $X$ is ${\cal T}(p)$-closed
if and only if $X=X^{p}$.

Proof. It is enough to show

(a) $\emptyset=\emptyset^{p}$;

(b) $\omega^{*}=(\omega^{*})^{p}$;

(c) $X_{0}\cup X_{1}=(X_{0}\cup X_{1})^{p}$ if
$X_{i}=X_{i}^{p}~(i=0,1)$ and

(d) $\cap_{i\in I}X_{i}=(\cap_{i\in I}X_{i})^{p}$ if each $X_{i}$
satisfies $X_{i}=X_{i}^{p}$.

Now (a) and (b) are obvious, as are the inclusions $\subseteq$ of (c)
and (d).

(c) $(\supseteq$)
If $f:\omega\approx f[\omega] \subseteq X_{0}\cup X_{1}$ satisfies
$\overline{f}(p)=x\in(X_{0}\cup X_{1})^{p}$ then with
$A_{i}=\{n<\omega:f(n)\in X_{i}\}$ we have $A_{0}\cup A_{1}\in p$
and hence $A_{\overline{i}}\in p$ for suitable $\overline{i}\in\{0,1\}$;
changing the values of $f$ on $\omega\backslash A_{\overline{i}}$
if necessary
(to ensure $f[\omega]\subseteq A_{\overline{i}}$),
we conclude that $x=\overline{f}(p)\in
X_{\overline{i}}^{p}=X_{\overline{i}}\subseteq X_{0}\cup X_{1}$.

(d) $(\supseteq$). If $x=\overline{f}(p)$ with
$f:\omega\approx f[\omega] \subseteq\cap_{i}X_{i}$ then
$x\in\cap_{i}(X_{i}^{p})=\cap_{i}X_{i}$.~~$\Box$

5.2. Remarks. (a) In the terminology of Lemma 5.1, the topology ${\cal
T}(p)$ is defined by the relation

${\cal T}(p)=\{\omega^{*}\backslash X:X\subseteq\omega^{*}, X$ is ${\cal
T}(p)$-closed\}.

(b) For notational convenience we denote by $I(p)$ the set of ${\cal
T}(p)$-isolated points of $\omega^{*}$, and we write
$A(p)=\omega^{*}\backslash I(p)$.
Clearly $ x\in I(p)$ if and only if $x$ is not a ``discrete limit"
of points in $\omega^{*}\backslash\{x\}$, that is, if and only if every
$f:\omega\approx f[\omega] \subseteq\omega^{*}\backslash\{x\}$ satisfies $\overline{f}(p)\neq
x$. The fact that $I(p)\neq\emptyset$ has been known for many years.
Indeed, Kunen~[10] has shown that there exist $2^{\bf c}$-many points
$x\in\omega^{*}$ such that $x\notin{\rm cl}_{\beta(\omega)}A$ whenever
$A\subseteq\omega^{*}\backslash\{x\}$ and $|A|\leq\omega$. (These are
the so-called weak-$P$-points of $\omega^{*}$.)

As a mnemonic device one may think of $A(p)$ and $I(p)$ as the sets
of $p$-accessible and $p$-inaccessible points, respectively.

(c) For $X\subseteq\omega^{*}$ the set $X^{p}$ may fail to be closed.
Indeed, the ${\cal T}(p)$-closure of $X\subseteq\omega^{*}$ is determined
by the following iterative procedure (cf. also [1]).

5.3. Lemma. Let $X\subseteq\omega^{*}$. For $\xi\leq\omega^{+}$ define
$X_{\xi}$ by :

$X_{0}=X$;

$X_{\xi}=\cup_{\eta<\xi}~X_{\eta}$ if $\xi$ is a limit ordinal;

$X_{\xi+1}=X_{\xi}^{p}$.\\
\noindent Then $X_{\omega^{+}}={\cal T}(p)-{\rm cl}~X$.~~$\Box$

The following fact, noted in [8], [5], [6], is crucial to many studies
of $\omega^{*}$ (see also [3](16.13) for a proof). One may
capture the thrust of this lemma by paraphrasing the
picturesque terminology of Frol\'{i}k~[6]: ``No type produces
itself."

5.4. Lemma. No homeomorphism from $\beta(\omega)$ into $\omega^{*}$ has
a fixed point.~~$\Box$

5.5. Lemma. Let $A$ and $B$ be countable, discrete subsets of
$\omega^{*}$, with $A\subseteq B^{*}$. Then $A^{p}\cap B^{p}=\emptyset$.

Proof. If $x\in A^{p}\cap B^{p}$ we may suppose without loss of
generality that there are $f:\omega\approx A$ and $g:\omega\approx
B$ such that $x=\overline{f}(p)=\overline{g}(p)$. The function
$h=f\circ g^{-1}:B\approx A\subseteq B^{*}$ satisfies

$\overline{f}\circ\overline{g^{-1}}=\overline{h}:\beta(B)\approx\beta(A)\subseteq
 B^{*}$\\
\noindent and $\overline{h}(x)=x\in B^{*}$, contrary to Lemma
5.4.~~$\Box$

5.6. Corollary. Let 
$A$ and $B$ be countably infinite, discrete subsets of $\omega^{*}$
such that $A\cap B=\emptyset$.
Then $A^{p}\cap B^{p}=\emptyset$.

Proof. Let $x\in A^{p}\cap B^{p}$ and let
$f:\omega\rightarrow f[\omega]\subseteq A$
and $g:\omega\rightarrow g[\omega]\subseteq B$
satisfy $x=\overline{f}(p)=\overline{g}(p)$.
Leaving $f$ and $g$ unchanged on suitably chosen elements of
$p$, but making modifications elsewhere if necessary, we assume
without loss of generality that either $f[\omega]\subseteq (g[\omega])^{*}$
or $g[\omega]\subseteq(f[\omega])^{*}$ or
$f[\omega]\cap(g[\omega])^{*}=(f[\omega])^{*}\cap g[\omega]=\emptyset$.
By Lemma 5.5 the first of these
possibilities, and by symmetry the second, cannot occur. We conclude that
$f[\omega]\cup g[\omega]$
is a countable, discrete subset of $\omega^{*}$ such that
$f[\omega]\cap g[\omega]=\emptyset$; it follows that
$(f[\omega])^{*}\cap (g[\omega])^{*}=\emptyset$, since
every countable (discrete) subset of $\omega^{*}$ is $C^{*}$-embedded (cf.
[7](14.27, 14N.5), [3](16.15).
This contradicts the relation
$x\in(f[\omega])^{*}\cap(g[\omega])^{*}$.~~$\Box$

5.7. Corollary. If $\omega^{*}\supseteq X\in{\cal T}(p)$, then
$X^{p}\in{\cal T}(p)$.

Proof. If $\omega^{*}\backslash X^{p}$ is not ${\cal T}(p)$-closed then
there is $f:\omega\approx f[\omega] =A\subseteq\omega^{*}\backslash X^{p}$ such that
$x=\overline{f}(p)\in X^{p}$. Since $X\in{\cal T}(p)$ we have $x\in
X^{p}\backslash X$ so there is $g:\omega\approx g[\omega]=B\subseteq
X$ such that $x=\overline{g}(p)$. From $A\cap B=\emptyset$
and 5.6 now follows
$x\in A^{p}\cap B^{p}=\emptyset$, a contradiction.~~$\Box$

5.8. Corollary. If $\omega^{*}\supseteq X\in{\cal T}(p)$ then ${\cal
T}(p)-{\rm cl}~X\in{\cal T}(p)$.

Proof. This is immediate from 5.3 and 5.7. $\Box$

Our goal is to 2-color the points of $\omega^{*}$ in such a way that
every copy of $\omega\cup\{p\}$ receives two colors. First we
consider how to extend a given coloring function.

5.9. Lemma. Let $\omega^{*}\supseteq X\in{\cal T}(p)$ and let
$c:X\rightarrow2=\{0,1\}$ be a function with no monochromatic copy of
$\omega\cup\{p\}$ (that is, if $X\supseteq Y\approx\omega\cup\{p\}$
then $c^{-1}(\{i\})\cap Y\neq\emptyset$ for $i\in\{0,1\}$). Then
$c$ extends to $\tilde{c}:X^{p}\rightarrow2$ with no monochromatic
copy of $\omega\cup\{p\}$.

Proof. Set $X_{i}=c^{-1}(\{i\})$ for $i\in2=\{0,1\}$, so that
$X^{p}=X_{0}^{p}\cup X_{1}^{p}$ by 5.1(c) and

$(X_{0}^{p}\backslash X)\cap(X_{1}^{p}\backslash X)=\emptyset$\\
\noindent by 5.6. Since $\{X~,~X_{0}^{p}\backslash X~,~X_{1}^{p}\backslash
X\}$ is a partition of $X$, the function
$\tilde{c}:X^{p}\rightarrow2$, given by the rule

$\tilde{c}(x)=c(x)$ if $x\in X$\\

\hspace{.35in}$=1$ if $x\in X_{0}^{p}\backslash X$

\hspace{.35in}$=0$ if $x\in X_{1}^{p}\backslash X$,\\
\noindent in well-defined. To see that $\tilde{c}$ is as required let
$h:\omega\cup\{p\}\approx A\cup\{x\}\subseteq X^{p}$ with
$h:\omega\approx A,~h(p)=x$. Modifying $h$ (as before) if necessary,
we assume without loss of generality that either (i)~$A\subseteq
X_{0}$ or (ii)~$A\subseteq X_{0}\backslash X$ (the cases $A\subseteq
X_{1}, A\subseteq X_{1}^{p}\backslash X$ are treated symmetrically).
In case~(i) we have $\tilde{c}\equiv0$ on $A$ and $\tilde{c}(x)=1$
(since either $x\in X$ or $x\in X_{0}^{p}\backslash X$); case~(ii)
cannot arise, since $x\in X$ violates $X\in{\cal T}(p)$ while $x\in
X^{p}\backslash X$ violates Corollary 5.6.~~$\Box$

Combining Lemmas 5.9 and 5.3 yields this.

5.10. Lemma. Let $\omega^{*}\supseteq X\in{\cal T}(p)$ and let
$c:X\rightarrow\{0,1\}$ be a function with no monochromatic copy of
$\omega\cup\{p\}$. Then $c$ extends to $\tilde{c}:{\cal T}(p)-{\rm
cl}~X\rightarrow\{0,1\}$ with no monochromatic copy of
$\omega\cup\{p\}$.~~$\Box$

The preceding lemma indicates how to extend a coloring function
from $X\in{\cal T}(p)$ over ${\cal T}(p)-{\rm cl}~X$, but it remains
to initiate the coloring procedure. For this purpose it is convenient to
consider a particular base ${\cal S}(p)$ for the topology ${\cal
T}(p)$. We call the elements of ${\cal S}(p)$ the $p$-{\it
satellite sets}.

5.11. Definition. Let $x\in\omega^{*}$. A set $S=S(x)$ is a $p$-{\it
satellite set} based at $x$ if there are a tree
$T\subseteq\omega^{<\omega}=\cup_{n<\omega}~\omega^{n}$ (ordered by
containment) and for $s\in T$ a point $x_{s}\in S$ and
$U_{s}\subseteq\omega^{*}$ such that

(i) $U_{s}$ is open-and-closed in the usual topology of $\omega^{*}$;

(ii) $x=x_{\langle\rangle}$ with
$\langle\rangle$ the empty sequence;

(iii) $U_{\langle\rangle}=\omega^{*}$;

(iv) if $x_{s}\in S(x)$ and $x_{s}\in A(p)$ then:
$\{x_{s\hat{~}n}:n<\omega\}$ enumerates the range of a function $f$
such that $f:\omega\approx f[\omega] \subseteq\omega^{*}$
with $\overline{f}(p)=x_{s}$, and
$\{U_{s\hat{~}n}:n<\omega\}$
is a pairwise disjoint family such that
$x_{s\hat{~}n}\in U_{s\hat{~}n}\subseteq U_{s}$;

(v) if $x_{s}\in S(x)$ and $x_{s}\in I(p)$ then $s$ is a maximal
node in $T$ (and $x_{s\hat{~}n},~U_{s\hat{~}n}$
are defined for no $n<\omega$).

5.12. Remark. It is not difficult to see that for every $x\in X\in{\cal
T}(p)$ there is $S=S(x)\in{\cal S}(p)$ such that $x\in S\subseteq
X$. (If $x\in I(p)$ one takes $S=\{x\}$; if $x_{s}\in S\cap X$
has been defined one uses (iv) and $X\in{\cal T}(p)$ to choose
$x_{s\hat{~}n}\in S\cap X$ if $x_{s}\in A(p)$.) That each of the sets
$S(x)$ is ${\cal T}(p)$-open is immediate from Corollary 5.6 above. It
follows that ${\cal S}(p)$ is indeed a base for  ${\cal T}(p)$.

5.13. Theorem. Every $p\in\omega^{*}$ satisfies $\omega^{*}\not\rightarrow(\omega\cup\{p\})^{1}_{2}$.

Proof. Let $\{S(x(i)):i\in I\}$ be a maximal pairwise disjoint
subfamily of ${\cal S}(p)$. For each $i\in I$ define
$c_{i}:S(x(i))\rightarrow2$ by

$c_{i}(x(i)_{s})=0$ if length of $s$ is even

\hspace{.6in}$=1$ if length of $s$ is odd.\\
\noindent It is clear from Corollary 5.6 that not only each function
$c_{i}$ on $S(x(i))$, but also the function

$c=\cup_{i\in I}~c_{i}:\cup_{i\in I}~S(x(i))\rightarrow2$,\\
\noindent is monochromatic on no copy of $\omega\cup\{p\}$. Since
$\cup_{i\in I}~S(x(i))$ is
${\cal T}(p)$-open and
${\cal T}(p)$-dense in $\omega^{*}$, the
desired result follows from Lemma 5.10.~~$\Box$

5.14. Theorem.  The relation $\omega^{*}\rightarrow(Y)^{1}_{2}$ fails for every infinite, countably
compact space $Y$.

Proof. Given infinite $Y\subseteq\omega^{*}$ there is $f:\omega\approx f[\omega] \subseteq Y$,
and if $Y$ is countably compact there is $p\in\omega^{*}$ such that
$\overline{f}(p)\in Y$. Since $f[\omega]$ is $C^{*}$-embedded in
$\omega^{*}$ we have

$\omega\cup\{p\}\approx f[\omega]\cup\{\overline{f}(p)\}\subseteq Y$,\\
\noindent so $\omega^{*}\not\rightarrow(Y)^{1}_{2}$ follows from
$\omega^{*}\not\rightarrow(\omega\cup\{p\})^{1}_{2}$.~~$\Box$

5.15. Remarks. (a) We cite three facts which (taken together) show that
the index set $I$ used in the proof of Theorem 5.13 satisfies
$|I|=2^{\bf c}$: (i)~ The set $W$ of weak-$P$-points of $\omega^{*}$
introduced by Kunen~[10]
satisfies $|W|=2^{\bf c}$; (ii)~each $S(x)\in{\cal S}(p)$ satisfies
$|S(x)|\leq\omega$; (iii)~$W\subseteq I(p)$, so
$W\subseteq\cup_{i\in I}~S(x(i))$.

(b) With no attempt at a complete topological classification, we note
five elementary properties enjoyed by each of our topologies ${\cal
T}(p)$ on $\omega^{*}$.

~~~~~(i) ${\cal T}(p)$ refines the usual topology of $\omega^{*}$, so
${\cal T}(p)$ is a Hausdorff topology.

~~~~~(ii) ${\cal T}(p)$ has $2^{\bf c}$-many isolated points. (Indeed, we
have noted already that the set $W$ of weak-$P$-points satisfies
$|W|=2^{\bf c}$ and $W\subseteq I(p)$.)

~~~~~(iii) Since ${\cal S}(p)$
is a base for ${\cal T}(p)$
and each $S(x)\in{\cal S}(p)$
satisfies $|S(x)|\leq\omega$, the topology ${\cal T}(p)$ is locally
countable.

~~~~~(iv) From Theorem 1.1(b) it is easy to see that if $S(x)\in{\cal S}(p)$
and $|S(x)|=\omega$, then $|{\cal T}(p)-{\rm cl}~S(x)|={\bf c}$.
Thus ${\cal T}(p)$ is not a regular topology for $\omega^{*}$.

~~~~~(v) According to Corollary 5.8, the ${\cal T}(p)$-closure of each ${\cal
T}(p)$-open subset of $\omega^{*}$ is itself ${\cal T}(p)$-open. Such
a topology is said to be extremally disconnected.

(c) In our development of ${\cal T}(p)$ and its properties we did not
introduce explicitly the Rudin-Frol\'{i}k pre-order $\sqsubseteq$
on $\omega^{*}$ (see [5], [6], or [13], or [3] for an expository
treatment) since doing so does not appear to simplify the arguments. We
note however (as in [1]) that the relation $\sqsubseteq$ lies close to
our work: For $x,p\in\omega^{*}$ one has $p\sqsubset x$ if and only if
some
$f:\omega\approx f[\omega] \subseteq\omega^{*}$ satisfies $\overline{f}(p)=x$.~~$\Box$

\S6. Questions.

Perhaps this paper is best viewed as establishing some boundary
conditions which may help lead to a solution of the following ambitious
general problem.

6.1. Problem. Characterize those spaces $Y$ such that $\omega^{*}\rightarrow(Y)^{1}_{2}$.~~$\Box$

There are $P$-spaces $Y$ such that $|Y|=2^{\bf c}$ and
$Y\subseteq_{h}\omega^{*}$. (For example, according to van Douwen's theorem cited
above, one may take $Y=P(\omega^{*})$.) According to Theorem 2.2, the
relation $\omega^{*}\rightarrow(Y)^{1}_{2}$ fails for each such $Y$. This situation suggests the
following question.

6.2. Question. Does $\omega^{*}\rightarrow(Y)^{1}_{2}$
for every $P$-space $Y$ such that $Y\subseteq_{h}\omega^{*}$ and
$|Y|<2^{\bf c}$? What if $|Y|={\bf c}$?~~$\Box$

We have no example of a non-$P$-space $Y$ such that $\omega^{*}\rightarrow(Y)^{1}_{2}$, so we are
compelled to ask:

6.3. Question. If $Y$ is a space such that $\omega^{*}\rightarrow(Y)^{1}_{2}$, must $Y$ be a
$P$-space?~~$\Box$

For $|Y|=\omega$, Question 6.3 takes this simple form:

6.4. Question. If $Y$ is a countable space such that $\omega^{*}\rightarrow(Y)^{1}_{2}$, must $Y$
be discrete?~~$\Box$

6.5. Remark. In connection with Question 6.4 it should be noted that
there exists a countable, dense-in-itself subset $C$ of $\omega^{*}$
such that every $x\in C$ satisfies

(*) $x\notin{\rm cl}_{\beta(\omega)}~D$ whenever $D$ is discrete and
$D\subseteq C\backslash\{x\}$\\
\noindent (equivalently: $\omega\cup\{p\}\subseteq_{h} C$
fails for every $p\in\omega^{*}$). To find such $C$ we follow the
construction of van Mill~[11](3.3, pp.~53-54). Let $E$ be the absolute
(i.e., the Gleason cover) of the Cantor set $\{0,1\}^{\omega}$, let
$\pi:E\rightarrow\{0,1\}^{\omega}$ be perfect and irreducible, and
embed $E$ into $\omega^{*}$ as a {\bf c}-OK set; then every countable
$F\subseteq\omega^{*}\backslash E$ satisfies $E\cap{\rm
cl}_{\beta(\omega)}F=\emptyset$. Now by the method of
[11](3.3) for $t\in\{0,1\}^{\omega}$ choose
$x_{t}\in\pi^{-1}(\{t\})$ such that every discrete
$D\subseteq E\backslash\{x_{t}\}$ satisfies $x_{t}\notin{\bf
cl}_{\beta(\omega)}D$, and take $C=\{x_{t}:t\in C_{0}\}$ with
$C_{0}$ a countable, dense subset of $\{0,1\}^{\omega}$. Since
$\pi$ is irreducible
the set $C$ is dense in $E$ and is dense-in-itself, and it is easy
to see that condition (*) is satisfied.

Of course no element of $C$ is a $P$-point of $\omega^{*}$. The
existence in ZFC of non-$P$-points $x\in\omega^{*}$ such that
$x\notin{\rm cl}_{\beta(\omega)}D$ whenever $D$ is a countable,
discrete, subspace of $\omega^{*}\backslash\{x\}$ is given explicitly by
van Mill~[11]; see also Kunen~[9] for a construction in ZFC~+~CH (or, in
ZFC~+~MA) of a set $C$ as above.

For the set $C$ constructed above the relation $\omega\cup\{p\}\subseteq_{h}C$
fails for every $p\in\omega^{*}$, so the following question, closely
related to Question 6.4, is apparently not answered by the methods of
this paper.

6.6. Question. Let $C$ be a countable, dense-in-itself subset of
$\omega^{*}$ such that $\omega\cup\{p\}\subseteq_{h} C$ fails for every
$p\in\omega^{*}$. Is the relation
$\omega^{*}\rightarrow(C)^{1}_{2}$ valid?~~$\Box$\\
\newpage
\begin{center}
{\large List of References}
\end{center}

[1] W. W. Comfort and Akio Kato. Non-homeomorphic disjoint spaces whose union
is $\omega^{*}$. Rocky Mountain J. Math. To appear.

[2] W. W. Comfort and Akio Kato. Topological partition relations of
the form $\omega^{*}\rightarrow(Y)^{1}_{2}$. Abstracts Amer. Math. Soc.
74 (1991), 288-289 (=~Abstract \#91T-54-25).

[3] W. W. Comfort and S. Negrepontis. The Theory of Ultrafilters.
Grundlehren der math. Wissenschaften Band 211. Springer-Verlag.
Berlin, Heidelberg and New York. 1974.

[4] Alan Dow and Jan van Mill. An extremally disconnected Dowker space.
Proc. Amer. Math. Soc. 86 (1982), 669-672.

[5] Zde\v{n}ek Frol\'{i}k.
Sums of ultrafilters. Bull. Amer. Math. Soc. 73 (1967), 87-91.

[6] Zde\v{n}ek Frol\'{i}k.
Fixed points of maps of $\beta{\bf N}$. Bull. Amer.
Math. Soc. 74 (1968), 187-191.

[7] Leonard Gillman and Meyer Jerison.
Rings of Continuous Functions. D. Van Nostrand Co., Inc. Princeton, Toronto,
New York and London. 1960.

[8] M. Kat\v{e}tov. A theorem on mappings. Commentationes Math. Universitatis
Carolinae 8 (1967), 431-433.

[9] Kenneth Kunen. Some points in $\beta{\bf N}$. Math. Proc.
Cambridge Philosophical Soc. 80 (1976), 385-398.

[10] Kenneth Kunen. Weak $P$-points in ${\bf N}^{*}$.
Colloquia Mathematica Societatis J\'{a}nos Bolyai 23 (1978), 741-749.

[11] Jan van Mill.
Sixteen topological types in $\beta\omega\backslash\omega$.
Topology and Its Applications 13 (1982), 43-57.

[12] Jan van Mill.
An introduction to $\beta\omega$.
In: Handbook of Set-theoretic Topology, pp. 503-567. Edited by K. Kunen and J.
E. Vaughan. North-Holland Publ. Co. Amsterdam, New York, Oxford. 1984.

[13] Mary Ellen Rudin.
Partial orders on the types in $\beta\bf N$.
Trans. Amer. Math. Soc. 155 (1971), 353-362.

[14] Stevo Todorcevic.
Partition Problems in Topology. Contemporary Mathematics, vol. 84.
American Mathematical Society. Providence, Rhode Island. 1989.

[15] William Weiss. Partitioning topological spaces. In: Mathematics
of Ramsey Theory, pp. 154-171. Edited by J. Ne\v{s}et\v{r}il and
V. R\"{o}dl. Springer-Verlag. Berlin. 1990.

[16] W. A. R. Weiss. Weiss's questions. In: Open Problems
in Topology, pp. 77-84. Edited by Jan van Mill and G. M. Reed.
North-Holland Publ. Co. Amsterdam. 1990.

\end{document}